\documentclass[a4paper,abstracton]{scrartcl}

\usepackage[utf8]{inputenc}

\usepackage{amsmath}
\usepackage{amssymb}

\usepackage{algpseudocode}
\usepackage{algorithm}

\usepackage[usenames,dvipsnames]{color}
\usepackage{graphicx}
\usepackage{subfigure}

\usepackage{cite}

\newcommand{\cU}{\mathcal{U}}
\newcommand{\cS}{\mathcal{S}}

\newcommand{\R}{\mathbb{R}}
\newcommand{\X}{\mathcal{X}}
\newcommand{\F}{\mathcal{F}}
\newcommand{\cA}{\mathcal{A}}

\newcommand{\SR}{{\rm S \hspace{-0.15em} R}}
\newcommand{\aSR}{{\rm a \hspace{-0.15em} S \hspace{-0.15em} R}}

\usepackage{authblk}
\usepackage[normalem]{ulem} 

\usepackage{framed}
\usepackage{rotating}

\definecolor{gray1}{gray}{.95}
\definecolor{gray2}{gray}{.5}
\newenvironment{Grau}{%
	\MakeFramed {\advance\hsize-\width \FrameRestore}}%
	{\endMakeFramed}

{
\begin{Grau}
\begin{center}
\begin{minipage}{0.95\textwidth}\centering\hspace{-1.2mm}\textsf{#1}\end{minipage}
\end{center}
\end{Grau}
}

\newtheorem{theorem}{Theorem}

\newtheorem{remark}[theorem]{Remark}

\begin{document}

\title{Algorithm Engineering in Robust Optimization\thanks{Partially supported by grant
SCHO 1140/3-2 within the DFG programme {\it Algorithm Engineering} and by FP7-PEOPLE-2009-IRSES under the
OptALI grant agreement no 246647.}}

\author[1]{Marc Goerigk\thanks{Email: goerigk@mathematik.uni-kl.de}}
\affil[1]{University of Kaiserslautern, Germany}

\author[2]{Anita Sch\"obel\thanks{Email: schoebel@math.uni-goettingen.de}}
\affil[2]{University of G\"ottingen, Germany}

\date{}

\maketitle

\begin{abstract}
Robust optimization is a young and emerging field of research having received
a considerable increase of interest over the last decade.
In this paper, we argue that the the algorithm engineering methodology fits 
very well to the field of robust optimization and
yields a rewarding new perspective on both the
current state of research and open research directions.

To this end we go through the algorithm engineering cycle of design and
analysis of concepts, development and implementation of algorithms,
and theoretical and experimental evaluation. 
We show that many ideas of algorithm engineering have already been applied in publications
on robust optimization. Most work on robust optimization is devoted to 
analysis of the concepts and the development of algorithms, some papers deal with the
evaluation of a particular concept in case studies, and work on
comparison of concepts just starts. 
What is still a drawback in many papers on robustness 
is the missing link to include the results of the
experiments again in the design. 
\end{abstract}

\section{Introduction}

Similar to the approach of {\it stochastic} optimization, robust optimization deals with models in which the
exact data is unknown, but bounded by a set of possible realizations (or scenarios). Contrary to 
the former approach, in robust optimization, one typically
refrains from assuming a given probability distribution over the scenarios. While the first steps in robust optimization
trace back to the work of Soyster \cite{Soyster}, it has not emerged as a field of research in its own right before the late 90s
with the seminal works of Ben-Tal, Nemirovski, and co-authors (see \cite{BenTalNemi1998,BenTalNemi1999}, and many more).

In this section, we first describe the general setting of robust optimization in more detail, and then discuss the
algorithm engineering methodology and its application, which gives 
a natural structure for the remainder of the paper.

\paragraph{Uncertain optimization problems.}

Nearly every optimization problem suffers from uncertainty to some
degree, even if this does not seem to be the case at first
sight. Generally speaking, we may distinguish two types of
uncertainty: {\it Microscopic} uncertainty, such as numerical errors
and measurement errors; and {\it macroscopic} uncertainty, such as
forecast errors, disturbances or other conditions changing the
environment where a solution is implemented.

In ``classic'' optimization, one would define a so-called {\it nominal
  scenario}, which describes the expected or ``most typical''
behavior of the uncertain data. Depending on the uncertainty type,
this scenario may be, e.g., the coefficient of the given precision for
numerical errors, the measured value for measurement errors, the most
likely forecast for forecast errors, or an average environment for
long-term solutions. Depending on the application, computing such a
scenario may be a non-trivial process, see, e.g.,
\cite{Jenkins2000275}.

In this paper we consider optimization problems that can be written in the form
\begin{align*}
(P) \hspace{1cm} \min\ &f(x) \\
\mbox{s.t. }   & F(x) \leq 0 \\
& x \in \X,
\end{align*}
where $F:\R^n \to \R^m$ describes the $m$ problem constraints, $f:\R^n \to \R$ is the objective function, and $\X \subseteq \R^n$ is the variable space. In real-world applications, both the constraints and the objective may depend on parameters which are uncertain. In order to accommodate such uncertainties, instead of $(P)$, the following parameterized {\it family} of problems is considered:
\begin{align*}
(P(\xi)) \hspace{1cm}\min\  &f(x,\xi) \\
\mbox{s.t. }   & F(x,\xi) \leq 0 \\
& x \in \X,
\end{align*}
where $F(\cdot, \xi):\R^n \to \R^m$ and $f(\cdot, \xi):\R^n \to \R$
for any fixed $\xi \in \R^M$. Every $\xi$ describes a \textit{scenario} that may occur.

\medskip
Although it is in practice often not known exactly which values such a scenario $\xi$ may take for an optimization problem $P(\xi)$,
we assume that it is known that $\xi$ lies within a
given \textit{uncertainty set} $\cU \subseteq \R^M$. Such an uncertainty set 
represents the scenarios which are likely enough to be considered.

The \textit{uncertain optimization problem}
corresponding to $P(\xi)$ is then
denoted as
\begin{equation}
\label{uncertain-prob}
 \left(P(\xi), \xi \in \cU\right).
\end{equation}

Note that the uncertain optimization problem in fact consists of a whole set of parameterized problems, that is often even infinitely large. The purpose of robust optimization concepts 
is to transform this family of problems back 
into a single problem, which is called the {\it robust counterpart}. 
The choice of the uncertainty set is of major impact not only for the respective application,
but also for
the computational complexity of the resulting robust counterpart.
It hence has to be chosen carefully by the modeler.

For a given uncertain optimization problem $(P(\xi), \xi \in \cU)$,
we denote by
\[ \F(\xi) =  \{x \in \X: F(x,\xi) \leq 0\}\]
the feasible set of scenario $\xi \in \cU$. Furthermore, if there exists a nominal scenario, it is denoted by $\hat{\xi}\in\cU$. \label{nominalscen} The optimal objective value for a single scenario $\xi\in\cU$ is denoted by $f^*(\xi)$. 

We say that an uncertain optimization problem $(P(\xi), \xi \in \cU)$ has 
convex (quasiconvex, affine, linear) uncertainty, when both functions, 
$F(x,\cdot):\cU \to \R^m$ and $f(x,\cdot):\cU \to \R$ are 
convex (quasiconvex, affine, linear) in $\xi$ for every fixed $x\in\X$.

\paragraph{Common uncertainty sets.}

There are some types of uncertainty sets that are frequently used in current literature. These include:
\begin{enumerate}
\item Finite uncertainty $\cU = \left\{\xi^1,\ldots,\xi^N\right\}$
\item Interval-based uncertainty $\cU = [\underline{\xi}_1,\overline{\xi}_1] \times \ldots \times [\underline{\xi}_M,\overline{\xi}_M]$
\item Polytopic uncertainty $\cU = \mathrm{conv}\left\{\xi^1,\ldots,\xi^N\right\}$
\item Norm-based uncertainty $\cU = \left\{ \xi \in\R^M : \Vert \xi - \hat{\xi} \Vert \le \alpha \right\}$ for some parameter~$\alpha \geq 0$
\item Ellipsoidal uncertainty $\cU = \left\{ \xi\in\R^M : \sqrt{\sum_{i=1}^M \xi^2_i/\sigma_i^2} \le \Omega\right\}$ for some parameter $\Omega\geq 0$
\item Constraint-wise uncertainty $\cU = \cU_1 \times \ldots \times \cU_m$, where $\cU_i$ only affects constraint $i$
\end{enumerate}
where ${\rm conv}\left\{\xi^1,\ldots,\xi^N\right\} = \left\{ \sum_{i=1}^N \lambda_i \xi^i : \sum_{i=1}^N \lambda_i =1, \lambda\in\R^N_+\right\}$ denotes the convex hull of a set of points.
Note that this classification is not exclusive, i.e., a given uncertainty set can belong to multiple types at the same time.

\paragraph{The algorithm engineering methodology, and the structure of this paper.}

\begin{figure}[htbp]
\centering\includegraphics[width=0.83\textwidth]{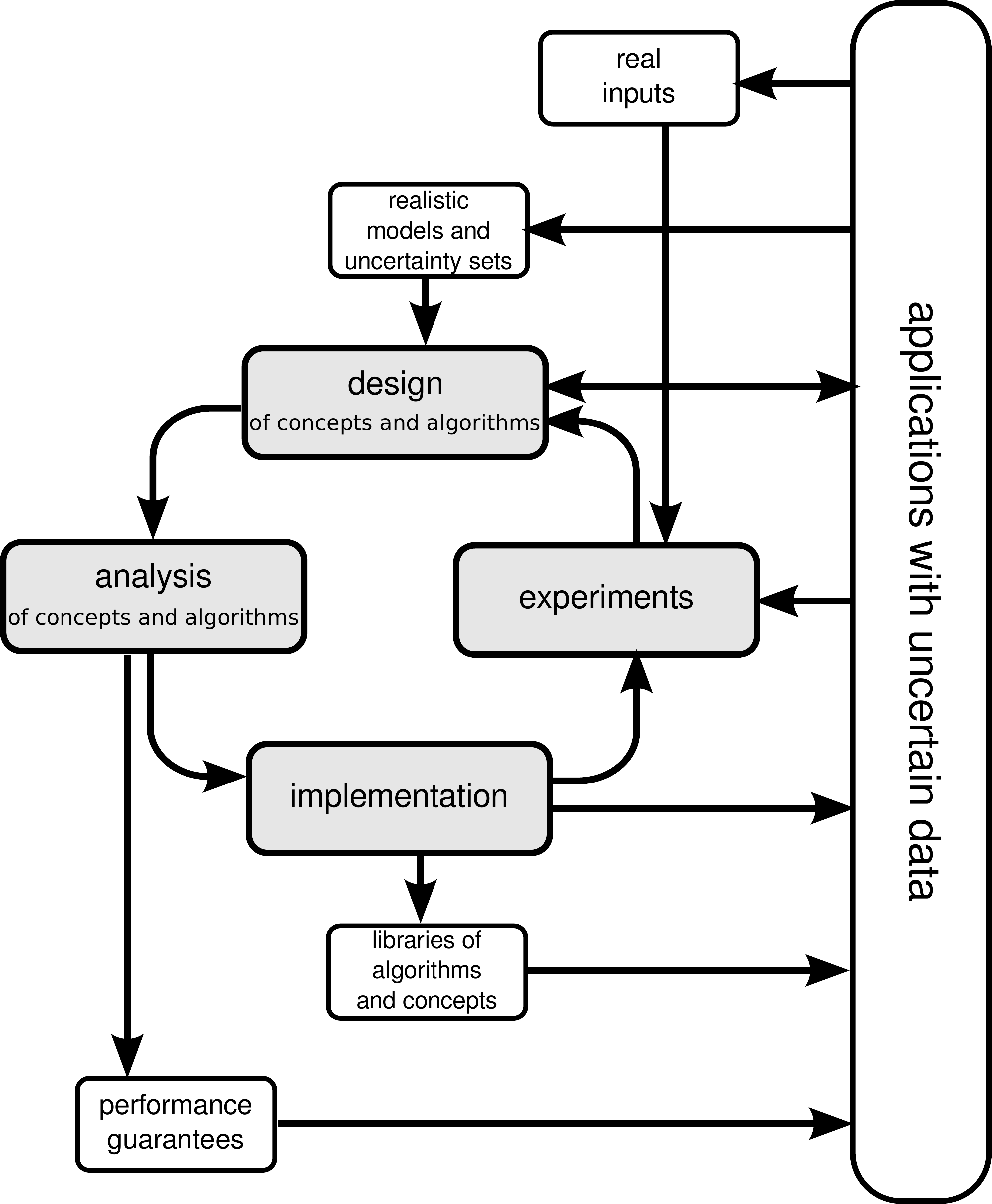}
\caption{The algorithm engineering cycle for robust optimization following \cite{San09}.}\label{ae-cycle}
\end{figure}

In the algorithm engineering approach, a feedback cycle 
between {\it design}, {\it analysis}, {\it implementations}, and {\it experiments} 
is used (see \cite{San09} for a detailed discussion). 
We reproduce this cycle for robust optimization in Figure~\ref{ae-cycle}.

While this approach usually focuses on the design and analysis of {\it algorithms}, one needs to consider the important role that different {\it concepts} play in robust optimization. 
Moreover, 
as is also discussed later, there is a thin line between what is to be considered a 
robustness concept, 
and an algorithm -- e.g.,  the usage of a simplified model for a robustness concept
could be considered as a new concept, but also as a heuristic algorithm for the original 
concept. We will therefore consider the design and analysis of both, concepts and algorithms.

The algorithm engineering approach has been successfully applied
to many problems and often achieved impressive speed-ups (as in routing algorithms, see, e.g. \cite{AEshortestpath} and the book \cite{muller2010algorithm}).

Even though this aspect has not been sufficiently acknowledged in the robust optimization 
community, the algorithm engineering paradigm fits very well in the line of research done 
in this area:
In algorithm engineering it is of particular importance that the single steps
in the depicted cycle are not considered individually, but that special structure
occurring in typical instances is identified and used in the development and analysis
of concepts and algorithms.
As we will show in the following sections these links to real-world applications and
to the structure of the uncertain data are of special importance in particular in
robust optimization. Various applications with different understandings of 
what defines a robust
solution triggered the development of the different robustness
concepts (see Section~\ref{sec:design}) while the particular structure of the uncertainty
set led to adapted algorithms (see Section~\ref{sec-reformulations}).

Moreover, the algorithm engineering cycle
is well-suited to detect the \emph{missing} research links to push the developed
methods further into practice. 
A key aspect of this paper hence is to draw further attention 
to the potential of algorithm engineering for robust optimization.

We structure the paper along the algorithm engineering cycle, where we discuss each step separately, 
providing a few exemplarily papers dealing with the respective matters. Missing links to 
trigger further research in this areas are pointed out. Specifically, we consider
\begin{itemize}
\item design of robustness concepts in Section~\ref{sec:design},
\item analysis of robustness concepts in Section~\ref{sec:analysis},
\item design and analysis of algorithms in Section~\ref{sec:algorithms}, and
\item implementations and experiments in Section~\ref{sec:experiments}.
\end{itemize}
Applications of robust optimization are various, and strongly 
influenced the design of robustness concepts while the design of algorithms was
rather driven by an analysis of the respective uncertainty sets. Some of these relations
are mentioned in the respective sections.
The paper is concluded in Section~\ref{sec:conclusion} where we also demonstrate
on some examples how the previously
mentioned results can be interpreted in the light of the algorithm engineering methodology.

\section{Design of Robustness Concepts}\label{sec:design}

Robust optimization started with rather conservative concepts hedging against everything
that is considered as being likely enough to happen. Driven by various other situations and applications calling for
``robust'' solutions these concepts were further developed. 
In this section we give an overview on the most important older and some recent concepts.
We put special emphasis on the impact applications with uncertain data
have on the design of robustness concepts (as depicted in in Figure~\ref{ae-cycle}), 
and how real-world requirements influence the development of robustness models.

\subsection{Strict Robustness}
\label{current-strict}

This approach, which is sometimes also known as classic robust optimization, one-stage robustness, min-max optimization, absolute deviation, or simply {\it robust optimization}, can be seen as the pivotal starting point in the field of robustness.
A solution $x\in\X$ to the uncertain problem $(P(\xi),\xi\in\cU)$
is called \textit{strictly robust} if it is feasible for all scenarios in $\cU$, i.e. if
$F(x,\xi) \leq 0$ for all $\xi \in \cU$. The objective usually follows the pessimistic view of minimizing the worst-case over all scenarios. Denoting the set of strictly robust solutions with respect to the uncertainty set $\cU$ by
\[\SR(\cU) = \bigcap_{\xi\in\cU} \F(\xi),\]
the strictly robust counterpart of the uncertain optimization problem is given as
\begin{eqnarray*}
\label{RC-1}
\mbox{(SR)} \hspace{1cm}
\min \ \ \sup_{\xi \in \cU} &f(x,\xi)\label{RC-11}\\
\mbox{s.t. } &x \in \SR(\cU)\\
&x\in\X.
\end{eqnarray*}
The first to consider this type of problems from the perspective of {\it generalized linear programs} was Soyster \cite{Soyster} for uncertainty sets $\cU$ of type  \[\cU=K_1 \times \ldots \times K_n,\]
where the set $K_i$ contains possible column vectors $A_i$ of the coefficient matrix~$A$. 

Subsequent works on this topic include \cite{falk1976} and \cite{thuente1980}.
\medskip

However, building this approach into a strong theoretic framework is due to a series of papers by Ben-Tal, Nemirovski, El Ghaoui and co-workers \cite{El-Ghaoui97,BenTalNemi1998,BenTalNemi1999,BenTalNemi2000}.
A summary of their results can be found in the book \cite{RObook}.
Their basic underlying idea is to hedge against all scenarios that may occur. As they argue, such an approach
makes sense in many settings, e.g., when constructing a bridge which must be stable, no matter which traffic
scenario occurs, or for airplanes or nuclear power plants. However, this high degree of 
conservatism of strict robustness is not applicable to all situations which call for robust solutions. 
An example for this is timetabling in public transportation: being strictly robust for a timetable means that all
announced arrival and departure times have to be met, no matter what happens. This may mean to add
high buffer times, depending on the uncertainty set used, and thus would not result in a practically applicable timetable. Such applications triggered
research in robust optimization on ways to relax the concept. 
We now describe some of these approaches.

\subsection{Cardinality Constrained Robustness}
\label{current-bs}

One possibility to overcome the conservatism of strict robustness is to shrink the uncertainty set $\cU$. This has been conceptually introduced by Bertsimas and Sim in \cite{BertSim04} for linear programming problems. Due to this reason, this concept is sometimes also known as ``the approach of Bertsimas and Sim'', sometimes also under the name ``$\Gamma$-robustness''. Analyzing the structure of uncertainty sets in typical applications, they observed
that it is unlikely that all coefficients
of one constraint change simultaneously to their worst-case values.
Instead they propose to hedge only against scenarios in which 
at most $\Gamma$ uncertain
parameters per constraint change to their worst-case values, i.e., they
restrict the number of coefficients which are allowed to change leading to
the concept of \emph{cardinality constrained robustness}.
Considering a constraint of the form
\[ a_1 x_1 + \ldots + a_n x_n \le b\]
with an uncertainty $\cU = \{a\in\R^n : a_i \in [\hat{a}_i -d_i, \hat{a}_i + d_i], i=1,\ldots,n\}$, 
their robustness concept requires a solution $x$ to satisfy
\[
\sum_{i=1}^n \hat{a}_i x_i + \max_{S\subseteq\{1,\ldots,n\},\atop{|S|=\Gamma}} \left\{ \sum_{i\in S} d_i|x_i|  \right\} \le b
\]
for a given parameter $\Gamma\in \{0,\ldots,n\}$.
Any solution $x$ to this model hence hedges against all scenarios in which at most 
$\Gamma$ many uncertain coefficients may deviate from their nominal values at the same time. 

It can be shown that
cardinality constrained robustness can also be considered as strict robustness using the convex hull of the 
\emph{cardinality-constrained uncertainty set}
\[ \cU(\Gamma)=\{a \in \cU:  a_i \not= \hat{a}_i \mbox{ for at most $\Gamma$ indices $i$} \} \subseteq \cU.\]

Since $conv(\cU(\Gamma))$ is a polyhedral set,
results on strict robustness with respect to polyhedral uncertainty can also be applied to 
cardinality constrained robustness.
\medskip

Note that this approach also extends to fractional values of $\Gamma$.
Their concept has been extended to uncertainty sets under general norms in \cite{Bertsimas2004510}. The approach to combinatorial optimization problems has been generalized in \cite{Atam06} and \cite{GoetzmannStillerTelha2011}.

\subsection{Adjustable Robustness}
\label{current-adjustable}

In \cite{BeGoGuNe2003} 
a completely different observation of instances occurring in real-world
problems with uncertain data is used: Often the variables can be decomposed
into two sets. The values for the \textit{here-and-now variables} have to be found by 
the robust optimization algorithm in advance, while the decision about the
\textit{wait-and-see variables} can wait until the actual scenario $\xi \in \cU$ becomes known.
Note that this is similar to two-stage programming in stochastic optimization.

We therefore assume that the variables $x=(u,v)$ are splitted into $u \in \X^1 \subseteq \R^{n_1}$ and $v \in \X^2 \subseteq \R^{n_2}$ with $n_1+n_2=n$, where the variables $u$ need to be determined before the scenario $\xi \in \cU$ becomes known, while the variables $v$ 
may be determined after $\xi$ has been realized. Thus, we may also write $x(\xi)$ to emphasize the dependence of $v$ on the scenarios. The uncertain optimization problem $(P(\xi),\xi\in\cU)$ is rewritten as
\begin{align*}
P(\xi)\quad \min\ & f(u,v,\xi)\\
&F(u,v,\xi)\le 0\\
&(u,v)\in \X^1\times \X^2.
\end{align*}
When fixing the here-and-now variables, one has to make sure that for any possible scenario $\xi \in \cU$ there exists $v \in \X^2$ 
such that $(u,v)$ is feasible for $\xi$.
The set of adjustable robust solutions is therefore given by
\begin{align*}
\aSR &= \{u\in \X^1 : \forall \xi\in\cU\ \exists v\in \X^2 \text{ s.t. } (u,v)\in\F(\xi)\}\\
&= \bigcap_{\xi\in\cU} Pr_{\X^1}(\F(\xi)),
\end{align*}
where $Pr_{\X^1}(\F(\xi))=\{u\in \X^1: \exists v\in \X^2 \text{ s.t. }(u,v) \in \F(\xi)\}$
denotes the projection of $\F(\xi)$ on $\X^1$.

The worst case objective for some $u \in \aSR$ is given as
\[z^{\aSR}(u) = \sup_{\xi\in\cU} \ \ \inf_{v: (u,v) \in \F(\xi)} f(u,v,\xi).\]
The adjustable robust counterpart is then given as
\[ \min \{z^{\aSR}(u): u \in \aSR\}. \]

Note that this setting is also useful for another type of
problem instances, namely, if auxiliary variables are used that do not represent decisions, e.g., additional
variables to model the absolute value of a variable.

There are several variations of the concept of adjustable robustness. Instead of two stages, multiple stages are possible. In the approach of finitely adaptable solutions \cite{5460988}, instead of computing a new solution for each scenario, a set of possible static solutions is computed, such that at least one of them is feasible in each stage.

Furthermore, the development of adjustable robustness was preceded by the similar approach of Mulvey et al \cite{mulvey1995}. They considered an uncertain linear optimization problem of the form
\begin{align*}
\mbox{(P(}B,C,e\mbox{))}\hspace{1cm} \min \ &c^tu + d^t v\\
\mbox{s.t. } &A u = b\\
&B u + C v = e\\
&u \in\R^{n_1}_+, v \in\R^{n_2}_+,
\end{align*}
where $u$ represents a vector of {\it design} variables that cannot be adjusted, and $v$ a vector of {\it control} variables that can be adjusted when the realized scenario becomes known. For a finite uncertainty set $\cU=\{(B^1,C^1,e^1),\ldots,(B^N,C^N,e^N)\}$, their robust counterpart is given as
\begin{align*}
\mbox{(Mul)}\hspace{1cm} \min \ &\sigma(u,v^1,\ldots,v^N) + \omega\rho(z^1,\ldots,z^N)\\
\mbox{s.t. } &A u = b\\
&B^i u + C^i v^i + z^i = e^i\ \forall i=1,\ldots,N\\
&u \in\R^{n_1}_+, v^i \in\R^{n_2}_+, z^i\in\R^m.
\end{align*}
The variables $z^i$ are introduced to measure the infeasibility in every scenario, i.e., the deviation 
from the right-hand side. The function $\sigma$ represents the {\it solution robustness}. It can be modeled as a worst-case function of the nominal objective
\[ \sigma(u,v^1,\ldots,v^N) = c^t u + \max_{i=1,\ldots,N} d^tv^i \]
or, when probabilities $p^i$ are known, as an expected nominal objective. The function $\rho$ on the other hand represents the {\it model robustness} and depends on the infeasibility of the uncertain constraints. Possible penalty functions are
\begin{align*}
\rho(z^1,\ldots,z^N) &= \sum_{i=1}^N p_i \sum_{j=1}^m \max\{0,z_j^i\}\\
\mbox{or } &= \sum_{i=1}^N p_i (z^i)^tz^i.
\end{align*}
As (Mul) is actually a bicriteria model, $\omega$ is used as a scalarization factor to combine both 
objectives.

\subsection{Light Robustness}
\label{current-light}

The lightly robust counterpart of an uncertain optimization problem, as developed in 
\cite{FischMona09} and generalized in \cite{Scho13} is again application driven. Originally developed
for timetabling, the idea of light robustness is that a solution must not be too bad in the 
nominal case. For example, the printed timetable should have short travel
times if everything runs smoothly and without disturbances; or a planned schedule should have
a small makespan. In this sense a certain nominal quality is fixed.
Among all solutions satisfying this standard, the concept asks for the most 
``reliable'' one with respect to constraint violation. Specifically, 
the general lightly robust counterpart as defined in \cite{Scho13} is of the following form:
\begin{align*}
\mbox{(LR)}\hspace{1cm}\min\ & \sum_{i=1}^m w_i \gamma_i\\
\mbox{s.t. } \ & f(x,\hat{\xi}) \le f^*(\hat{\xi}) + \rho\\
& F(x,\xi) \le \gamma \qquad \forall \xi \in \cU\\
&x \in \X, \gamma \in\R^m,
\end{align*}
where $w_i$ models a penalty weight for the violation of constraint $i$ and $\rho$ determines 
the required nominal quality. We denote by $\hat{\xi}$ the nominal scenario, as introduced on page~\pageref{nominalscen}.
This approach was in its first application in \cite{FischMona09} used as a further development of the
concept of cardinality constrained robustness (see Section~\ref{current-bs}).

Note that a constraint of the form $F(x,\xi)\le 0$ is equivalent to a constraint $\lambda F(x,\xi) \le 0$ for any $\lambda > 0$; therefore, the coefficients $w_i$ play an important role in balancing the allowed violation of the given constraints.

\subsection{Recoverable Robustness}
\label{current-recovery}

Similar to adjustable robustness, {\it recoverable robustness} is again a two-stage concept.
It has been developed in \cite{CDDFN07,sebastiandiss,LLMS09,DDN09}
and has independently also been used in \cite{Savel09}. Its basic idea is to allow a class of \textit{recovery algorithms} $\mathcal A$ that can be used in case of a disturbance. A solution $x$ is called \textit{recovery robust} with respect to $\mathcal A$ if for any possible scenario $\xi \in \cU$
there exists an algorithm $A \in \mathcal A$ such that $A$ applied to the solution $x$ and the
scenario $\xi$ constructs a solution $A(x,\xi) \in \F(\xi)$, i.e., a solution which is feasible for the
current scenario.

The recovery robust counterpart according to \cite{LLMS09} is the following:
\begin{align*}
\mbox{(RR)}\hspace{1cm} \min_{(x,A)\in\F(\hat\xi)\times\cA} &f(x)\\
\mbox{s.t. }\ &A(x,\xi) \in \F(\xi)\ \forall \xi\in\cU.
\end{align*}
It can be extended by including the recovery costs of a solution $x$: Let $d(A(x,\xi))$ be a possible vector-valued function that measures the costs of the recovery, and let $\lambda\in\Lambda$ be a limit on the recovery costs, i.e., $\lambda\ge d(A(x,\xi))$ for all $\xi\in\cU$. Assume that there is some cost function $g: \Lambda \to \R$ associated with $\lambda$.

Setting
\[A(x,\xi,\lambda) \in \F'(\xi) \iff d(A(x,\xi)) \le \lambda\ \wedge\ A(x,\xi)\in\F(\xi)\]
gives the recovery robust counterpart with limited recovery costs:
\begin{align*}
\mbox{(RR-LIM)}\hspace{1cm} \min_{(x,A,\lambda)\in\F(\hat\xi)\times\cA\times\Lambda} &f(x) + g(\lambda)\\
\mbox{s.t. }\ &A(x,\xi,\lambda) \in \F'(\xi)\ \forall \xi\in\cU.
\end{align*}

Due to the generality of this robustness concept, the computational tractability heavily depends on the problem, the recovery algorithms and the uncertainty under consideration. In 
\cite{atmos2010,Tapas2011-Goerigk-Sch,marcdiss,GoeSchoe13a},
the concept of recoverable robustness has been considered under the usage of metrics to 
measure recovery costs. 
The aim is to minimize the costs when recovering, where they differ between recovering 
to a feasible
solution (``recovery-to-feasibility''), and recovering to an optimal solution (``recovery-to-optimality'') 
in the realized scenario.

\subsection{Regret Robustness} 
\label{current-regret}

The concept of regret robustness differs from the other presented robustness concepts insofar it 
usually only considers uncertainty in the objective function. Instead of minimizing the worst-case performance of a solution, it minimizes the difference to the objective function of the best solution that would have been possible in a scenario.
In some publications, it is also called \emph{deviation robustness}.

Let $f^*(\xi)$ denote the best objective value in scenario $\xi\in\cU$. 
The min-max regret counterpart of an uncertain optimization problem with uncertainty in 
the objective is then given by
\begin{align*}
\mbox{(Regret)}\hspace{1cm} \min\ \sup_{\xi\in\cU}\ &\Big( f(x,\xi) - f^*(\xi) \Big)\\
\mbox{s.t. }\ &F(x) \le 0\\
&x\in\X.
\end{align*}
Regret robustness is a concept with a vast amount of applications, e.g., in location theory or in scheduling.
For a survey on this concept, see \cite{Aissi2009} and \cite{KouYu97}. In a similar spirit, the concept of lexicographic $\alpha$-robustness has been recently proposed \cite{Kalai2012722}. Its basic idea is to evaluate a fixed solution by reordering the set of scenarios according to the performance of the solution. This performance curve is then compared to an ideal curve, where the optimization problem is solved separately for
every scenario.

\subsection{Some Further Robustness Concepts}
\label{current-further}

\paragraph{Reliability.}

Another approach to robust optimization is to relax the constraints of strict robustness. This leads to the concept of
\textit{reliability}
\label{reliable}
of Ben-Tal and Nemirovski \cite{BenTalNemi2000}, in which
the constraints $F(x,\xi) \leq 0$ are replaced by
$F(x,\xi) \leq \gamma$ for some $\gamma \in \R^m_{\geq 0}$. A solution $x$
which satisfies
\[ F(x,\xi) \leq \gamma \mbox{ for all } \xi \in \cU \]
is called \textit{reliable with respect to $\gamma$}. The goal is to find a reliable solution
which minimizes the original objective function in the worst case.
Similar to light robustness, one has to be careful that the representation of the constraints does not affect the reliability of the solution, otherwise
one may obtain the counter-intuitive result that, although the constraints 
$F(x,\xi) \leq 0$ can also be written as $\Psi(F(x,\xi)) \leq 0$ for any increasing $\Psi$ with $\Psi(0)=0,$
what is understood by a robust solution may be different if one models the constraints with $F$ or with $\Psi(F)$.

\paragraph{Soft Robustness.}
The basic idea of {\it soft robustness} as introduced in \cite{soft} is to handle the conservatism of the strict robust approach by considering a nested family of uncertainty sets, and allowing more deviation 
in the constraints for larger uncertainties. Specifically, instead of an uncertainty set $\cU\subseteq\R^M$, a family of uncertainties $\{\cU(\varepsilon)\subseteq\cU\}_{\varepsilon>0}$ with $\cU(\varepsilon_1) \subseteq \cU(\varepsilon_2)$ for all $\varepsilon_2 \ge \varepsilon_1$ is used. The set of soft robust solutions is then given as
\[{\rm soft\hspace{-0.15em} R} = \left\{ x\in\X : F(x,\xi) \le \varepsilon \ \forall \xi\in\cU(\varepsilon),\ \varepsilon >0\right\}.\]
Note that strict robustness is a special case with $\cU(\varepsilon)=\cU$ for all $\varepsilon>0$.

In \cite{soft}, the authors show that a solution to the softly robust counterpart  -- i.e., the optimization over ${\rm soft\hspace{-0.15em} R}$ with a worst-case objective -- can be found by solving a sequence of strictly robust counterparts
using a bisection approach over $\varepsilon$, and analyze the numerical performance on a bond portfolio 
and an asset allocation problem.

\paragraph{Comprehensive Robustness.}

While the adjustable robust approach relaxes the assumption that all
decisions have to be made before the realized scenario becomes known,
the approach of comprehensively robust counterparts \cite{BeBoNe06} also
removes
the assumption that only scenarios defined in the uncertainty
set $\cU$ need to be considered. Instead, using a distance measure $dist$ in
the space of scenarios, and a distance measure $\overline{dist}$ in the solution space,
they assume that the further away a scenario is from the uncertainty
set, the further away the corresponding solution is allowed to be from
the set of feasible solutions. As in adjustable robustness, the dependence of $x$ on the scenario $\xi$ is allowed, 
and we may write $x(\xi)$. The adjustable robust counterpart
is extended to the following problem:
\begin{align*}
\mbox{(CRC)}\hspace{1cm} \min \ &z\\
\mbox{s.t. } &f(x(\xi),\xi) \le z + \alpha_0 dist(\xi,\cU) \ \forall \xi \\
&\overline{dist}(x(\xi),\F(\xi)) \le \alpha dist(\xi,\cU) \ \forall \xi,
\end{align*}
where $\alpha,\alpha_0$ denote sensitivity parameters. This formulation needs further formal specification, which is provided in \cite{BeBoNe06}.

\paragraph{Uncertainty Feature Optimization.}

Instead of assuming that an explicit uncertainty set is given, which may be hard to model for real-world problems, the uncertainty feature optimization (UFO) approach \cite{ufo} rather assumes that the robustness of a solution is given by an explicit function. For an uncertain optimization problem $(P(\xi))$, let $\mu: \R^n \to \R^p$ be a measure for $p$ robustness features. The UFO-counterpart of the uncertain problem is then given by
\begin{align*}
\mbox{(UFO)}\hspace{1cm} {\rm vecmax}\  &\mu(x)\\
\mbox{s.t. } &F(x) \le 0\\
&f(x) \le (1+\rho) f^*(\hat{\xi})\\
&x\in\X,
\end{align*}
where $f^*(\hat{\xi})$ denotes the best objective value to the nominal problem. The authors 
show that this approach generalizes both stochastic optimization and the concept of 
cardinality constrained robustness
of Bertsimas and Sim.

\subsection{Summary}
\label{unified}

As this section shows, we cannot actually speak of one concept 
or point-of-view to be ``robust optimization''; instead, we should see it as a vast collection of different 
robustness concepts, 
each providing their unique advantages and disadvantages. Generally speaking, there is usually a trade-off between the degree of freedom a concept gives to react to disruptions (including what is considered as being a disruption, i.e., the choice of the uncertainty  set), and its computational complexity. From an algorithm engineering point of view, the size of this ``toolbox'' of different concepts significantly helps with finding 
a suitable robustness concept
for a given problem. However, as these concepts are usually application-driven, they 
lack a generalizing systematics: Applications tend to develop ``their own approach'' to robustness instead of making use of the existing body of literature, and develop their own notation and names along the way. In fact, the very same concepts are known under plenty of names. Summaries as \cite{RObook,Aissi2009,bertsimas-survey,Roy2010629} usually avoid this Babylonian ``zoo'' of robustness concepts and nomenclature by focusing only on the mainstream concepts. Thus, we suggest the following pointer to further research:

\begin{remark}
Robust optimization needs a unified classification scheme.
\end{remark}

\section{Analysis of Robustness Concepts}\label{sec:analysis}

Not only the development of robustness concepts, but also their analysis is
data-driven. This becomes in particular clear when looking at the structure of the
underlying uncertainty set. A large amount of research in the analysis of robustness
concepts is devoted to finding equivalent 
problem formulations that are better tractable, using the structure of
the uncertainty set.

In this section we first review this line of research, and then briefly point out
exemplarily which other types of structure or ideas have been used in the analysis
of concepts.

\subsection{Using the structure of the uncertainty set}
\label{sec-reformulations}

\paragraph{Finite uncertainty set.}
If the uncertainty set $\cU=\{\xi^1,\ldots,\xi^N\}$ is a finite set containing not too many scenarios, 
most of the robustness concepts can be formulated as mathematical programs by just
adding the constraints for each of the scenarios explicitly. This can straightforwardly
been done for strict robustness yielding 
\begin{eqnarray*}
\mbox{(SR)} \hspace{1cm}
\min \ \ z \\
\mbox{s.t. } f(x,\xi^i) & \leq & z \ \mbox{ for } i=1,\ldots,N \\
\mbox{s.t. } F(x,\xi^i) & \leq & 0 \ \mbox{ for } i=1,\ldots,N \\
&& x  \in\X.
\end{eqnarray*}
as the strictly robust counterpart. Reliability and light robustness can be treated analogously.
In all three cases, the robust counterpart keeps many properties of the original (non-robust)
problem formulation: If the original formulation was e.g., a linear program, also its robust
counterpart is. The same holds for differentiability, convexity, and many other properties.

For regret robustness one needs to precompute the best objective function value
for each scenario $\xi^1,i=1,\ldots,N$ in order to receive again a straightforward reformulation.
Also in adjustable and recoverable robustness mathematical programming formulations
can be derived by adding a wait and see variable, or a group of 
recovery variables for \emph{each} of the scenarios.
This usually leads to a high number of additional variables but is (at least for linear programming)
often still solvable. 

Note that the concept of 
cardinality constrained robustness
does not make much sense for a finite set of scenarios since
it concerns the restriction \emph{which} scenarios might occur. For a finite set, scenarios
in which too many parameters change can be removed beforehand. 

\paragraph{Polytopic uncertainty.}

Let $f(x,\cdot)$ and $F(x,\cdot)$ be quasiconvex in $\xi$ for any fixed $x \in \X$. 
Then there are robustness concepts in which the following \emph{reduction result}
holds: The robust counterpart w.r.t. an uncertainty
set $\cU'$ is equivalent to the robust counterpart w.r.t. $\cU:=conv(\cU')$. 
In such cases the robust counterpart w.r.t. a polytopic uncertainty set
 $\cU=conv\{\xi^1,\ldots,\xi^N\}$ is equivalent to the robust counterpart w.r.t. 
the finite uncertainty set $\{\xi^1,\ldots,\xi^N\}$, hence 
the formulations for finite uncertainty sets can be used to treat polytopic uncertainties. 

We now review for which robustness concepts the reduction result holds.
First of all, this is true for strict robustness, 
For affine and convex uncertainty this was
mentioned in \cite{BenTalNemi1998}; the generalization to quasiconvex uncertainty is straightforward.
One of the direct consequences, namely that 
the robust counterpart of an uncertain 
linear program under these conditions is again a linear program was mentioned
in \cite{BenTalNemi2000}. The same result holds for reliability since the reliable robust counterpart
can be transformed to a strictly convex counterpart by defining 
$\tilde{F}(x,\xi)=F(x,\xi) - \gamma$.
For light robustness, the result is also true, see \cite{Scho13}. For the case of adjustable robustness,
\cite{BeGoGuNe2003} 
showed that the result holds for problems with fixed recourse. Otherwise, 
counterexamples can be constructed. 
The generalization to nonlinear
two-stage problems and quasiconvex uncertainty is due to \cite{takeda08}.
For recoverable robustness there exist special cases in which the
recovery robust counterpart is equivalent to an adjustable robust counterpart with fixed recourse. In
these cases, the result of \cite{BeGoGuNe2003} may be applied. However,
in general, recoverable robustness does
not allow this property. This also holds for recovery-to-optimality.

\paragraph{Interval-based uncertainty.}
Interval-based uncertainty can be interpreted as a special case of polytopic uncertainty
where the polytope is a box 
$\cU = [\underline{\xi}_1,\overline{\xi}_1] \times \ldots \times [\underline{\xi}_M,\overline{\xi}_M]$
with $2^M$ extreme points $(\xi_1,\xi_2,\ldots,\xi_M)^t \in \R^M$, where $\xi_i \in \{\underline{\xi}_i, \overline{\xi}_i\}, i=1,\ldots,M$.
Hence, all the results for polytopic uncertainty apply. They can often be simplified
by observing that not all extreme points are needed since the respective constraints 
often dominate each other, yielding a drastic speed-up when solving the robust 
counterpart.

For their concept of 
cardinality constrained robustness,
Bertsimas and Sim \cite{BertSim04}
considered interval-based uncertainty sets for linear programs. This
can be interpreted as strict robustness with a new uncertainty set $\cU'$ only allowing
scenarios in which not more than $\Gamma$ uncertain parameters per constraint change their values 
(see also \cite{Bertsimas2004510}).
This uncertainty set $\cU'$ is a polytope, hence the robust counterpart for 
cardinality constrained robustness
stays a linear program for interval-based uncertainty.

\paragraph{Ellipsoidal uncertainty.}
The case of ellipsoidal uncertainty is studied extensively for strict robustness 
and for adjustable robustness in \cite{RObook}.
It could be shown that often the constraint
\[ F(x,\xi) \leq 0 \mbox{ for all } \xi \in \cU \]
can be replaced by a finite number of constraints for ellipsoidal uncertainty sets. However,
it has been shown in \cite{RObook} that for ellipsoidal uncertainty,
the structure of the strictly robust counterpart gets more complicated. For example (see
\cite{BenTalNemi1998}) the strictly robust counterpart of a linear program 
is a conic quadratic program, the strictly robust counterpart of a quadratic constrained 
quadratic program is a semidefinite program, the strictly robust counterpart of a second order 
cone program is a semidefinite program, and the strictly robust counterpart of a semidefinite 
program is NP-hard. As mentioned before, all these results can be transferred to
reliability. 

For light robustness, it has been shown in \cite{Scho13} that the lightly robust counterpart of a linear
program with ellipsoidal uncertainty becomes a quadratic program. Ellipsoidal uncertainty could receive
more attention also for other robustness concepts (e.g., for regret robustness, which 
usually only considers finite or interval-based uncertainty, 
see \cite{Aissi2009}), or for adjustable robustness, see \cite{RObook}.

\subsection{Using duality}
\label{sec-dual}

Duality in uncertain programs has been considered as early as 1980,
see \cite{thuente1980}.  In \cite{Beck20091}, it is shown that ``the
primal worst equals the dual best'', i.e., under quite general
constraints, the dual of a strictly robust counterpart (a min-max problem)
amounts to optimization under the best case instead of the worst-case (a max-min problem).
Since
then, duality in robust optimization has been a vivid field of
research, see, e.g., \cite{Jeyakumar2013331} and \cite{Suzuki2013257}.
In the following, we highlight two applications of duality for robust optimization: One
for constraints, and one for objectives.

\subsubsection*{Duality in the constraints.}

Duality is a useful tool for the reformulation of robust constraints. We exemplarily
demonstrate this using two applications.

In \cite{BertSim04},
the authors show that the 
cardinality constrained robust
counterpart
can be linearized by using the dual of the inner maximization problem. This yields
\begin{align*}
&\sum_{i=1}^n \hat{a}_i x_i + z\Gamma + \sum_{i=1}^n p_i \le b\\
&z + p_i \ge d_iy_i \qquad \forall i=1,\ldots,n\\
&-y_i \le x_i \le y_i \qquad \forall i=1,\ldots,n\\
&p, y, z \ge 0.
\end{align*}

Note that a general, robust constraint of the form
\[ f(x,\xi) \le 0 \ \forall \xi\in\cU \]
can be rewritten as
\[ \max_{\xi\in\cU} f(x,\xi) \le 0.\]
This is used in \cite{nonlinear12}. With a concave function $f(x,\cdot)$ and an uncertainty set $\cU = \{ \hat{\xi} + A\zeta : \zeta\in Z\}$ with a nonempty, convex and compact set $Z$, applying duality yields
\[ \hat{\xi}^t v + \delta^*(A^Tv|Z) - f_*(v,x) \le 0 \]
where $\delta^*$ is the support function, $f_*$ is a conjugate function, and other technical requirements are met. This gives a very general tool to compute robust counterparts; e.g., a linear constraint of the form $f(x,\xi) = \xi^tx - \beta$ and $Z= \{\zeta : \Vert \zeta \Vert_2 \le \rho \}$ yields the counterpart $\hat{\xi}^t x + \rho \Vert A^tx\Vert_2 \le \beta$.

\subsubsection*{Duality in the objective.}

In many papers, duality is used to change the typical min-max objective of a robust counterpart
into a min min objective by using the dual formulation of the inner maximization problem. 

This method was first applied to the spanning tree problem \cite{Yaman200131}, and later extended to the general case of 
optimization problems with zero duality gap in \cite{Aissi2009}. Let an uncertain optimization problem of the form
\begin{align*}
\min\ &c^tx \\
\text{s.t. } x\in \X &= \{ x\in\{0,1\}^n : Ax \ge b \}
\end{align*}
with interval-based uncertainty in $c$ be given; i.e., $c_i \in [ \underline{c}_i , \overline{c}_i]$. Then we may write
\begin{align*}
&\min_{x\in\X} \max_{c\in\cU} \left( f(x,c) - f^*(c) \right) \\
&= \min_{x\in\X} \max_{c\in\cU, y\in\X} \left( c^tx - c^ty \right)\\
&= \min_{x\in\X} \left( \overline{c} x - \min_{y\in\X} c^{wc}(x) y \right),
\end{align*}
where $c^{wc}(x)$ denotes the regret worst-case for $x$, given as $\overline{c}_i$ if $x_i = 1$, and $\underline{c}_i$ if $x_i =0$. Using that the duality gap is zero, we can insert the dual to the inner optimization problem, and get the following equivalent problem:
\begin{align*}
\min\ &\overline{c}x - b^t y\\
\text{s.t. } & Ax \ge b \\
&A^ty \le (\overline{c} - \underline{c}) x + \underline{c}\\
&x\in \{0,1\}^n, y\in \R^n_+
\end{align*}
This reformulation can then be solved using, e.g., a branch and bound approach.

\section{Design and Analysis of Algorithms}\label{sec:algorithms}

Concerning the design and analysis of algorithms we concentrate on the most mature 
concept, namely on algorithms for strict robustness. Many approaches, often based on
similar ideas, also exist for regret optimization -- e.g., cutting plane approaches \cite{Inuiguchi1995526,Mausser1999157,ITOR:ITOR389}, or preprocessing considerations \cite{Yaman200131,Kasperski2010680}.
For the other concepts, approaches are currently still being developed.

The robust counterpart per se is a semi-infinite program; thus, all methods that apply to semi-infinite programming \cite{Lopez2007491} can be used here as well. However, the special 
min-max structure of the robust counterpart
allows improved algorithms over the general case, in particular for the reformulations 
based on special uncertainty sets
as mentioned in Section~\ref{sec-reformulations}.

In the following, we discuss algorithms that are generically applicable to strictly robust optimization problems.

\subsection{Finite Scenarios}
\label{1-RC-finite}

The case we consider here is that $\cU=\{\xi^1,\ldots,\xi^N\}$ is a finite set; i.e., the strictly robust counterpart (SR)
can be rewritten as
\begin{align*}
\min \max_{i=1,\ldots,N} &f(x,\xi^i)\\
\text{s.t. } &F(x,\xi^i) \le 0 \ \forall i=1,\ldots,N\\
&x \in \X
\end{align*}
Due to the lack of structure in the uncertainty set, these instances can be hard so solve, even though they have a similar structure as the nominal problem. 

\subsubsection{Branch and bound using surrogate relaxation.}

The following approach was introduced by \cite{KouYu97} for discrete optimization problems with uncertainty only 
in the objective: For any vector $\mu\in\R_+^N$, the surrogate relaxation SRC($\mu$) of (SR)
with uncertain objective function is given by
\begin{align*}
\mbox{SRC($\mu$) \hspace{1cm}}  
\min\ & \frac{1}{\sum_{\xi\in\cU} \mu_\xi} \sum_{\xi\in\cU} \mu_\xi f(x,\xi) \\
\text{s.t. } &F(x) \le 0 \\
&x \in \X
\end{align*}

Note that the structure of the nominal problem is preserved, which allows the usage of 
specialized
algorithms already known. Furthermore, the optimal objective value SRC$^*$($\mu$) of this problem is a lower bound on the optimal objective value 
$\SR^*$ of the strictly robust counterpart; 
and as the set of feasible solutions is the same, also an upper bound is provided by solving SRC($\mu$).

This approach is further extended by solving the problem
\[ \max_{\mu\in\R_+^N}  \text{SRC}^*\text{(}\mu\text{)}, \]
i.e., by finding the multiplier $\mu$ that yields the strongest lower bound. This can be done using a sub-gradient method.

The lower and upper bounds generated by the surrogate relaxation are then used within a branch and bound framework on the $x$ variables. The approach was further improved for the knapsack problem in \cite{Iida99,Taniguchi98}.

\subsubsection{Local search heuristics.}

In \cite{Sbihi2010339}, a local search-based algorithm for the knapsack problem with uncertain objective function is developed. 
We briefly list 
the main aspects. It makes use of two different search procedures: Given a feasible solution $x$ and a list of 
local neighborhood moves $M$, 
let $GS(x,M)$ (the generalized search) determine the worst-case objective value of every move, and return the best move along with its objective value. Furthermore, let $RS(x,M,S)$ (the restricted search) perform 
a random search using the moves $M$ with at most $S$ steps.

The \textit{cooperative local search} algorithm (CLS) works as follows: It first constructs
a heuristic starting solution, e.g., by a greedy approach. In every iteration, a set of moves $M$ 
is constructed using either the generalized search for sets with small cardinality, 
or the restricted search for sets with large cardinality. 
When a maximum number of iterations is reached, the best feasible solution found so far is returned.

\subsubsection{Approximation algorithms.}

A discussion of approximation algorithms for strict robustness
with finitely many scenarios is given, e.g., in \cite{Aissi2007281}, 
where it is shown that there is an 
FPTAS for the shortest path, the spanning tree, and the knapsack problem when the number of scenarios is constant; but the shortest path problem is not $(2-\epsilon)$-approximable, the spanning tree problem is not $(\frac{3}{2}-\epsilon)$-approximable, and the knapsack problem is not approximable at all when the number of scenarios is considered as a non-constant input.

The basic idea for their results is to use the relationship between the strictly robust counterpart (SR)
and multi-objective optimization: At least one optimal solution for (SR) is also an efficient solution in the multi-objective problem where every scenario is an objective. Thus, if the multi-objective problem has a polynomial-time $\alpha$-approximation algorithm, then also (SR) has a polynomial-time $\alpha$-approximation.

There exist many more approximation algorithms for specific problems.
For example, in \cite{feige07}, robust set covering problems are considered with implicitly given, exponentially many scenarios. For a similar setting of exponentially many, implicitly given scenarios for robust network design problems (e.g., Steiner tree), \cite{mirrokni13} presents approximation results. 
Approximation results using finite scenario sets for two-stage robust covering problems, 
min-cut and shortest path can be found in
\cite{ravi05} and \cite{ravi14}.

\subsection{Infinite Scenarios}

\subsubsection{Sampling.}

When we cannot make use of the structure of $\cU$ (i.e., 
when the reformulation approaches from Section~\ref{sec:analysis} cannot be applied,
or when we do not have a closed description of the set available),
we can still solve (SR) heuristically using a finite subset of scenarios (given that we have some sampling method available). The resulting problem can 
be solved using the algorithms described in Section~\ref{1-RC-finite}.

In a series of paper \cite{random1,random2,random3,random4}, the probability of a solution calculated by a sampled scenario subset being feasible for all scenarios is considered.
It is shown
that for a convex uncertain optimization problem, the probability of the 
\emph{violation event} $V(x)=P\{\xi\in\cU: x \notin \F(\xi)\}$ can be bounded by
\[P(V(x^*) > \epsilon) \le
\sum_{i=0}^{n-1} \begin{pmatrix}N\\i\end{pmatrix} \epsilon^i
(1-\epsilon)^{N-i},\] where $N$ is the sample size, $x^*\in\R^n$ is
an optimal solution with respect to the sampled
scenarios, and $n$ is (as before) the dimension of the
decision space. Note that the left-hand side is the probability of a
probability; this is due to fact that $V(x)$ is a random variable in
the sampled scenarios. In other words, if a desired probability of
infeasibility $\epsilon$ is given, the accordingly required sample size
can be determined. This result holds under the assumption that every
subset of scenarios is feasible, and is independent of the probability
distribution which is used for sampling over $\cU$.

As the number of scenarios sampled this way may be large, the sequential optimization approach \cite{FuWa07,FuWa09b,FuWa09} uses sampled scenarios one by one. Using the above probability estimates, a solution generated by this method is feasible for (SR) only within a certain probability. 
The basic idea is the following: 
We consider the set $\cS(\gamma)$ of feasible solutions with respect to a given quality level $\gamma$, i.e.,
\begin{align*}
\cS(\gamma) &= \{ x\in\X : f(x) \le \gamma, F(x,\xi) \le 0 \ \forall \xi\in\cU \} \\
&= \{ x \in \X : \nu(\gamma,x,\xi) \le 0 \ \forall \xi\in\cU \}
\end{align*}
where
\[ \nu(\gamma,x,\xi) = \left( \max \{ 0, f(x) - \gamma \}^2 + \max \{ 0, F(x,\xi)\}^2 \right)^{1/2} \]
Using a subgradient on $\nu$, the current solution is updated in every iteration using the sampled scenario $\xi$. Lower bounds on the number of required iterations are given to reach a desired level of solution quality and probability of feasibility.

\subsubsection{Outer-approximation and cutting-plane methods.}

For this type of algorithm, the general idea is to iteratively a) solve 
a robust optimization problem with a finite subset of scenarios, and b) 
use a worst-case oracle that optimizes over the uncertainty set $\cU$ for a given 
solution $x$. These steps can be done either exactly or approximately. 

Algorithms of this type have often been used, see, e.g., \cite{Reemtsen199487,MuBo09,BuNoAl13,bender11,robustbus,Montemanni20061479,fimo2012}; sometimes even without knowledge that such an approach already 
exists (see also the lacking unification in robust optimization mentioned in
Section~\ref{unified}).

The following general results should be mentioned.
\cite{MuBo09} show that this method converges under certain assumptions, and present further variations that improve the numerical performance of the algorithm. Cutting-plane methods are compared to compact formulations on general problem benchmarks 
in~\cite{fimo2012}. 
In \cite{BuNoAl13}, the implementation is considered in more detail: A distributed algorithm version is presented, 
in which each processor starts with a single uncertain 
constraint, and generated cutting planes are communicated.

\subsection{Algorithms for Specific Problems}

The goal of this section is to show how much one can benefit by using the structure a specific problem might have.
To this end, we exemplarily chose three specialized algorithms: The first solves an NP-hard problem in pseudo-polynomial time,
the second is a heuristic for another NP-hard problem,  and the third is a polynomial-time solution approach.
Note that many more such algorithms 
have been developed.

In \cite{Monaci20132625}, a dynamic programming algorithm is developed for the robust knapsack problem with 
cardinality constrained uncertainty
in the weights. Extending the classic dynamic programming scheme to also include the number of items that are on their upper bounds, they are able to show a $\mathcal{O}(\Gamma n c)$ time complexity, where $n$ is the number of items, and $c$ is the knapsack budget (note that this is not a polynomial algorithm). 
The key idea of the dynamic program is an easy feasibility check of a solution, which is achieved by using an item sorting based on the upper weight bound $\bar{w}_i$.
In computational experiments, instances with up to 5000 items can be solved in reasonable time.

The problem of min-max regret shortest paths with interval uncertainty is considered in \cite{Montemanni20041667}. The general idea is based on path ranking, and the conjecture that a path that ranks good on the worst-case scenario, may also rank good with respect to regret. Considering paths with respect to their worst-case performance order, they formulate a stopping
criterion when the regret of a path may not improve anymore. Note that the regret of a single path can in this case
easily be computed by assuming the worst-case length for all edges in the path, and the best-case length for all other edges. 
Experiments show a strong correlation between computation times and length of the optimal path.

While the former two problems are NP-hard (for regret shortest path, see \cite{Zielinski2004570}), a polynomial-time algorithm for the min-max regret 1-center on a tree with uncertain edge lengths and node weights is presented in \cite{Averbakh2000292}. A 1-center is a point on any edge of the tree for which the maximal weighted distance to all nodes is minimized. The algorithm runs in $\mathcal{O}(n^6)$ time, which can be reduced to $\mathcal{O}(n^2\log(n))$ for the unweighted case. It is based on the observation that an edge that contains an optimal solution can be found in $\mathcal{O}(n^2\log(n))$ time; however, determining its exact location for the weighted case is more complicated.

Further algorithms to be mentioned here are the polynomial algorithm for min-max regret flow-shop scheduling with two jobs from \cite{Averbakh2006761}; the polynomial algorithm for the min-max regret location-allocation problem from \cite{Conde20071025}; the heuristic for regret spanning arborescences from \cite{Conde2007561}; the polynomial algorithm for the min-max regret gradual covering location problem from \cite{Berman2011233}; and the PTAS for two-machine flow shop scheduling with discrete scenarios from \cite{Kasperski201236}.

\subsection{Performance Guarantees}

We now discuss performance guarantees in robust optimization. Measuring the performance of a robust solution or algorithm 
can be either done by developing
guarantees regarding the performance of an algorithm or of a heuristic solution; but 
also by developing performance guarantees that compare the solutions generated by different 
robustness concepts.

On the algorithmic side, standard measures like the approximation ratio (i.e., the ratio between the robust objective value of the heuristic and the optimal robust solution) can be applied. 
There are simple, yet very general approximation algorithms presented in \cite{Aissi2009} for strict robustness
and regret robustness: If the original problem is polynomially solvable, there is an $N$-approximation algorithm for finite uncertainty sets, where $N$ is the number of scenarios. Furthermore, there is a $2$-approximation algorithm for regret robustness with interval-based uncertainty \cite{Kasperski2006} by using the mid-point scenario. These results have been extended in \cite{Conde2012452}, see 
also the approximability survey \cite{Aissi2007281} on 
strict and regret robustness. 
We do not know of approximation algorithms for other robustness concepts, which would provide 
interesting insight in the structural differences between the robust counterparts.

Regarding the comparison between solutions generated by different concepts,
an interesting approach is to consider the quality of a strictly robust solution when used in an adjustable setting, as done in 
\cite{bertsimas2010geometric,bertsimas2011geometric}. The authors are able to develop performance guarantees solely based on the degree of symmetry of the uncertainty set.

Concerning the evaluation of a robust solution (and not the algorithm to compute it), there is no general consent how to proceed, and surprisingly little systematic research can be found regarding this field. The so-called \emph{robustness gap} as considered in \cite{BenTalNemi1998} is defined as the difference between the worst-case objective of the robust solution, and the worst optimal objective value over all scenarios, i.e., as 
$\SR^* - \sup_{\xi\in\cU} f^*(\xi)$, where $\SR^*$ 
denotes the optimal value of (SR). They show that in the case of constraint-wise affine uncertainty, a compact set $\X$, and some technical assumptions, this gap equals zero. However, the most widely used approach is computing the so-called {\it price of robustness} \cite{BertSim04}, which is usually defined as the ratio between the robust solution value and the nominal solution value, i.e., as
\[ \frac{\min_{x\in\SR} \sup_{\xi\in\cU} f(x,\xi)}{\min_{x\in\F(\hat{\xi})} f(x,\hat{\xi})} \]
As an example, \cite{MP13} presents the analytical calculation of the price of robustness for knapsack problems. 
Using an interval-based uncertainty set on the weights (i.e., the weight of item $i$ is in $[w_i - \underline{w}_i, w_i + \overline{w}_i]$) and a cardinality constrained robustness approach, they show that the price of robustness equals $1/(1+\lceil \delta_{max}\rceil)$ for $\delta_{max} := \max_i \overline{w}_i/w_i$ and $\Gamma = 1$. For $\Gamma \ge 2$, the price of robustness becomes $1/(1+\lceil 2\delta_{max}\rceil)$.

Note that this is a rather pessimistic view on robustness, as it only concentrates on the additional {\it costs} 
of a robust solution compared to the nominal objective function value of an optimal solution for the nominal case.
However, if the application under consideration is affected by uncertainty, the nominal solution will not necessarily find nominal conditions, hence
the robust solution may actually {\it save} costs compared to the nominal solution (which easily may be even infeasible). There is no general
``golden rule'' that would provide a fair evaluation for the performance of a robust solution.

Note that such a bound is not the kind of performance guarantee that was actually considered in \cite{BertSim04}. Rather, they developed probability bounds for the feasibility of a solution to the cardinality constrained approach depending on $\Gamma$. Using such bounds they argue that the nominal performance of a solution can be considerably increased 
without decreasing the probability of being feasible too much.

Summarizing the above remarks, we claim that:

\begin{remark}
Performance guarantees are not sufficiently researched in 
robust optimization.
\end{remark}

\section{Implementation and Experiments}\label{sec:experiments}

\subsection{Libraries}

In the following, we present some libraries that are designed for robust optimization. A related overview can also be found in \cite{ropipaper}.

\paragraph{AIMMS for Robust Optimization.}

AIMMS \cite{aimms}, which stands for ``Advanced Interactive Multidimensional Modeling System'', is a proprietary software that contains an algebraic modeling language (AML) for optimization problems.
AIMMS supports most well-known solvers, including Cplex\footnote{http://www-03.ibm.com/software/products/en/ibmilogcpleoptistud}, Xpress\footnote{http://www.fico.com/en/products/fico-xpress-optimization-suite} and Gurobi\footnote{http://www.gurobi.com/}.

Since 2010, AIMMS has offered a robust optimization add-on, which was developed in a partnership with A. Ben-Tal. The extension only considers the concepts of strict and adjustable robustness as introduced in Sections~\ref{current-strict} 
and \ref{current-adjustable}. 
As uncertainty sets, interval-based uncertainty sets, polytopic uncertainty sets, or ellipsoidal uncertainty sets
are supported and transformed to mathematical programs as described in Section~\ref{sec-reformulations}.
The respective transformations are automatically done when the model is translated from the algebraic modeling language to the solver.

\paragraph{ROME.}

While AIMMS focuses on the work of Ben-Tal and co-workers, ROME \cite{rome} (``Robust Optimization Made Easy'') takes its origins in the work of Bertsimas, Sim and co-workers. ROME is built in the MATLAB\footnote{http://www.mathworks.com/products/matlab/} environment, which makes it on the one hand intuitive to use for MATLAB-users, but on the other hand lacks the versatility of an AML. As a research project, ROME is free to use. It currently supports Cplex, Xpress and SDPT3\footnote{http://www.math.nus.edu.sg/{\textasciitilde}mattohkc/sdpt3.html} as solver engines.

ROME considers polytopic and ellipsoidal uncertainty sets, that can be further specified using the mean support, the covariance matrix, or directional deviations. Assuming an affine dependence of the wait-and-see variables, it then transforms 
the uncertain optimization problem to an adjustable robust counterpart. The strictly robust counterpart is included
as a special case.

\paragraph{YALMIP.}

Similar to ROME, YALMIP \cite{yalmip} is a layer between MATLAB and a solver that allows the modeling of optimization problems under uncertainty. Nearly all well-known solvers are supported, including Cplex, Gurobi and Xpress.

YALMIP considers strict robustness. In order to obtain the strict robust counterpart of an uncertain optimization
problems so-called {\it filters} are used: When presented a model with uncertainty, the software checks if one of these filters applies to generate the strictly
robust counterpart. Currently, five of these automatic transformations are implemented. A duality filter (which adds dual variables according to Section~\ref{sec-dual}), 
an enumeration filter for finite and polytopic scenario sets (which simply lists all relvant constraints), an explicit maximization filter (where a worst-case scenario is used), the P\'olya filter (which 
is based on an inner approximation of the set of feasible solutions), and an elimination filter (which sets variables affected by uncertainty to 0 and is used as a last resort).

\paragraph{ROPI.}

The Robust Optimization Programming Interface (ROPI) \cite{ropipaper,ropihp} is a C++ library that provides wrapper MIP 
classes to support a range of solvers. Using these generic classes, a robust counterpart is automatically generated 
given the desired robustness concept and uncertainty set. Contrary to the previous libraries, 
a wider choice of robustness concepts is provided: These include strict robustness,  adjustable robustness, light robustness, and different versions of recoverable robustness.
\bigskip

Even though a user can pick and choose between multiple robust optimization libraries, there is to the best of our knowledge no library of robust optimization {\it algorithms} available. All of the above implementations are based on reformulation approaches, which makes it possible to draw upon existing solvers. However, as described in Section~\ref{sec:algorithms}, 
there are plenty of specifically designed algorithms for robust optimization available. Making them readily-implemented 
available to the user should be a significant concern for future work in robust optimization.

\begin{remark}
There is no robust optimization library available with specifically designed algorithms other than reformulation approaches.
\end{remark}

\subsection{Applications}

As already stated, robust optimization has been application-driven; thus, there are abundant papers dealing with applications of some robustness approach to real-world or at least realistic problems. Presenting an exhaustive list would go far beyond the scope of this paper; examples include circuit design \cite{mani06}, emergency logistics \cite{BenTal20111177}, 
and load planning \cite{loadplanning} for adjustable robustness; supply chain optimization \cite{bert2006} and furniture planning \cite{JoseAlem2012139} for 
cardinality constrained
robustness; inventory control for comprehensive robustness \cite{Aharon2009922}; timetabling \cite{FischMona09,fischetti09}, and timetable information \cite{BGMHSS11} for light robustness; shunting \cite{CDDFN07}, timetabling \cite{robustness-overview08,atmos2010},
and railway rolling stock planning \cite{Kroon10} for recoverable robustness; and airline scheduling for UFO \cite{ufo-diss}.

Hence, we can state:

\begin{remark}
Robust optimization is application-driven.
\end{remark}

\subsection{Comparative Experiments}

In this section we consider research that either compares two robustness concepts 
to the same problem, 
or two algorithms for the same problem and robustness concept. We present a list of papers on the former aspect 
in Table~\ref{exp-table}, and a list of papers on the latter aspect in Table~\ref{exp-table2}. We do not claim completeness 
for these tables; rather, they should be considered as giving a general impression on recent directions of research.

\begin{sidewaystable}
\centering
\begin{tabular}{r|rr|lp{6.2cm}l|lp{6cm}}
Year & Paper &\hspace{2mm} & \hspace{2mm} & Problem  & \hspace{2mm} & \hspace{2mm} & Robustness Concept\\[1mm]
\hline&&&&&&\\[-3mm]
2008 & \cite{Bertsimas20083} &&& Portfolio management &&& strict and cc\\[2mm]
2009 & \cite{Aharon2009922} &&& Inventory control &&& adjustable and comprehensive \\[2mm]
2009 & \cite{Yin2009470} &&& Road improvement &&& strict and scenario-based \\[2mm]
2010 & \cite{atmos2010} &&& Timetabling &&& strict, buffered, light, and variations of recoverable\\[2mm]
2010 & \cite{KazemiZanjani2010882} &&& Sawmill planning &&& Mulvey with different recourse costs \\[2mm]
2010 & \cite{Xu2010707} &&& Water sensor placement &&& strict and regret \\[2mm]
2011 & \cite{Tapas2011-Goerigk-Sch} &&& LP Benchmarks &&& strict and recoverable \\[2mm]
2011 & \cite{Adasme20111377} &&& Wireless network resource allocation &&& finite and interval-based \\[2mm]
2011 & \cite{Lin2011361} &&& Newsvendor &&& strict and regret \\[2mm]
2013 & \cite{rotti} &&& Timetable information &&& strict and light\\[2mm]
2013 & \cite{rotti2} &&& Timetable information &&& strict and recoverable\\[2mm]
2013 & \cite{loadplanning} &&& Load planning &&& strict and adjustable\\ [2mm]
2013 & \cite{Agra2013856} &&& Vehicle routing &&& strict and adjustable \\[2mm]
\end{tabular}
\caption{Papers presenting experiments comparing at least two different 
robustness concepts. ``cc'' abbreviates ``cardinality constrained''.}\label{exp-table}
\end{sidewaystable}

\begin{sidewaystable}
\centering
\begin{tabular}{r|rr|lp{3.6cm}l|lp{2cm}l|lp{6cm}}
Year & Paper &\hspace{2mm} & \hspace{2mm} & Problem  & \hspace{2mm} & \hspace{2mm} & Concept  & \hspace{2mm} & \hspace{2mm} & Algorithms\\[1mm]
\hline&&&&&&\\[-3mm]
2005 & \cite{Montemanni2005771} &&& Spanning tree &&& regret &&& Branch and bound, MIP \\[2mm]
2006 & \cite{Montemanni20061479} &&& Spanning tree &&& regret &&& Bender's decomp., MIP, branch and bound \\[2mm]
2008 & \cite{nikulin08} &&& Spanning tree &&& regret &&& Simulated annealing, branch and bound, Bender's decomp.  \\[2mm]
2008 & \cite{Taniguchi98} &&& Knapsack &&& strict &&& branch and bound with and without preprocessing \\[2mm]
2008 & \cite{GonzalezVelarde2008797} &&& Capacitated sourcing &&& adjustable &&& tabu search \\[2mm]
2009 & \cite{Conde2009235} &&& Critical path &&& regret &&& MIP and heuristic \\[2mm]
2010 & \cite{deFariasJR20101610} &&& Machine scheduling &&& strict &&& MIP with and without cuts \\[2mm]
2010 & \cite{Bohle2010245} &&& Wine harvesting &&& cc robust &&& MIP and scenario generation \\[2mm]
2010 & \cite{Ng2010557} &&& Lot allocation &&& strict &&& branch-and-price and heuristics \\[2mm]
2011 & \cite{Catanzaro20111610} &&& Shortest path &&& regret &&& IP with and without preprocessing \\[2mm]
2011 & \cite{Pereira20111153} &&& Assignment &&& regret &&& MIP, Bender's decomp., genetic algorithms \\[2mm]
2012 & \cite{minmax-tabu} &&& Spanning tree &&& regret &&& tabu search and IP \\[2mm]
2012 & \cite{fimo2012} &&& diverse &&& cc robust &&& MIP and cutting planes \\[2mm]
2012 & \cite{Song20121988} &&& Knapsack &&& strict &&& local search and branch and bound \\[2mm]
2013 & \cite{Monaci20132625} &&& Knapsack &&& cc robust &&& dynamic programming and IP \\[2mm]
2013 & \cite{Ouorou2013318} &&& Capacity assignment &&& adjustable &&& approximations \\[2mm]
\end{tabular}
\caption{Papers presenting experiments comparing at least two algorithms for the same robustness concept. ``cc'' abbreviates ``cardinality constrained''.}\label{exp-table2}
\end{sidewaystable}

We conclude the following from these tables and the accompanying literature: 
Firstly, papers considering real-world applications that compare different robustness concepts are relatively rare. 
Applied studies are too often satisfied with considering only one approach of the many that are possible. 
Secondly, algorithmic comparisons dominantly stem from the field of min-max regret, where at the same time mostly
academic problems are considered. The efficient calculation of solutions for other robustness concepts is still a relatively 
open and promising field of research. 
Summarizing, we claim that:

\begin{remark}
There are too few comparative studies in robust optimization.
\end{remark}

A different aspect Table~\ref{exp-table} reveals 
is that most computational studies comparing at least two robustness concepts include strict robustness as a ``baseline concept''; accordingly, and unsurprisingly, the more tailor-made approaches will show an improved behavior for the application at hand. This is much similar to frequently published papers on optimization problems which compare a problem-specific method to a generic MIP solver, 
usually observing a better performance of the former compared to the latter. 

However, while a standard MIP solver
is often still competitive to problem-tailored algorithms, a robustness concept which does not capture the problem specifics at hand will nearly always be the second choice to one which uses the full problem potential.

\subsection{Limits of Solvability}

We show the approximate size of benchmark instances used for testing exact algorithms for a choice of robust problems in Table~\ref{exp-table3}. These values should rather be considered as rough indicators on the current limits of solvability than the exact limits themselves, as problem complexities are determined by many more aspects\footnote{Number of items for finite, strict knapsack is estimated with the pegging test from \cite{Taniguchi98}.}.

\begin{table}[htbp]
\centering
\begin{tabular}{rr|rr|ll|ll}
Problem && Approach && Size && Source\\[1mm]
\hline&&&&&&\\[-3mm]
Spanning tree  && interval regret && $\sim100$ nodes && \cite{Perez14}\\[1mm]
Knapsack && finite strict && $\sim 1500$ items && \cite{G14}\\[1mm]
Knapsack && \ finite recoverable && $\sim500$ items && \cite{BKK11}\\[1mm]
Knapsack && cc strict && $\sim5000$ items && \cite{Monaci20132625}\\[1mm]
Knapsack && cc recoverable && $\sim200$ items && \cite{BKK11b}\\[1mm]
Shortest path && interval regret && $\sim1500$ nodes && \cite{ChasseinGoerigk2014}\\[1mm]
Assignment && interval regret && $\sim500$ items && \cite{Pereira20111153}\\[1mm]
\end{tabular}
\caption{Currently considered problem sizes for exact algorithms.}\label{exp-table3}
\end{table}

What becomes immediately obvious is that these limits are much smaller than for their nominal problem counterparts, which can go easily into the millions.

\subsection{Learning from Experiments}

We exemplarily show how experimental results can be used to design better algorithms for robust optimization; thus, we highlight the potential that lies in following the algorithm engineering cycle. To this end, we consider the regret shortest path problem: Given a set of scenarios consisting of arc lengths in a graph, find a path from a fixed source node to a fixed sink node which minimizes the worst-case length difference to an optimal path for each scenario.

From a theoretical perspective, the problem complexity is well-understood. For discrete uncertainty sets (and already for only two scenarios), the problem was shown to be NP-hard in the seminal monograph \cite{KouYu97}. For interval-based uncertainty, \cite{Zielinski2004570} showed its NP-hardness.

Furthermore, it is known that the regret shortest path problem with a finite, but unbounded set of scenarios is not approximable within $2-\epsilon$. For the interval-case, a very simple 2-approximation algorithm (see \cite{Kasperski2006}) is known: All one needs to do is to compute the shortest path with respect to the midpoint scenario, i.e., the arc lengths which are the midpoint of the respective intervals.

To solve the interval regret problem exactly, a branch-and-bound method has been proposed \cite{Montemanni20041667}, which branches along the worst-case path in the graph. However, computational experience shows that the midpoint solution -- despite being ``only'' a 2-approximation -- is already an optimal, or close-to-optimal solution for many of the randomly generated benchmark instances.

Examining this aspect in more detail, \cite{ChasseinGoerigk2014} developed an instance-dependent approximation guarantee for the midpoint solution, which is always less or equal to 2, but usually lies around $\sim1.6-1.7$.

Using these two ingredients -- the strong observed performance of the midpoint solution, and its instance-dependent lower bound -- the branch-and-bound algorithm of \cite{Montemanni20041667} can be easily adapted, by using a midpoint-path-based branching strategy instead of the worst-case path, and by using the improved guarantee as a lower bound. The resulting algorithm considerably outperforms the previous version, with computation times two orders of magnitude better for some instance classes.

These modifications were possible by studying experimental results, improving thereupon the theoretical analysis, and feeding this analysis back to an algorithm. It is an example for the successful traversal of an algorithm engineering cycle, and we believe that many more such algorithmic improvements can be achieved this way.

\section{Algorithm Engineering in Robust Optimization and Conclusion}\label{sec:conclusion}

In this paper we propose to use the algorithm engineering methodology to better understand the open problems and challenges in robust optimization. Doing so, we were able to point out links between algorithm engineering and robust optimization, and we presented an overview on the state-of-the-art from this perspective.

In order to further stress the usefulness of the
algorithm engineering methodology, we finally 
present three examples. Each of them is composed of a series of
papers, which together follow the algorithm engineering cycle in robust optimization.

\emph{Example 1:} Development of new models based on shortcomings of previous ones.

\cite{Soyster} introduced the concept of strict robustness. 
This concept was illustrated in several examples (e.g. from linear programming, see \cite{BenTalNemi2000}
or for a cantilever arm as in \cite{BenTalNemi1998}) and analyzed for these examples
in a mathematical way. The analysis in these papers showed that the problem complexity
increases then introducing robustness (e.g., 
the robust counterpart of an uncertain linear program with ellipsoidal uncertainty
is an explicit conic quadratic program). Moreover, the authors recognized that
the concept is rather conservative introducing an approximate robust counterpart with a more moderate
level of conservatism. These ideas were taken up \cite{BertSim04} to start the next run through the algorithm
engineering cycle by introducing their new concept of cardinality constrained robustness, 
which is less conservative and computationally better tractable, but may be applied only to easier uncertainty sets. 
Applying this concept to train 
timetabling and performing experiments with it was the starting point of \cite{FischMona09} who
relaxed the constraints further and developed the concept of light robustness which was then later generalized
to arbitrary uncertainty sets by \cite{Scho13}.

\medskip

\emph{Example 2:} From one-stage to two-stage robustness.

Recognizing that the concept of strict robustness 
is too conservative, \cite{BeGoGuNe2003} proposed the first two-stage robustness approach by introducing
their concept of adjustable robustness. When applying this concept to several application of railway
planning within the ARRIVAL project (see \cite{ARRIVAL}), 
\cite{LLMS09} noted that the actions allowed to adjust a timetable
do not fit the practical needs. This motivated them to integrate recovery actions in robust planning yielding
the concept of recoverable robustness. Unfortunately, recovery robust solutions are hard to obtain. Research on
developing practical algorithms is still ongoing. recent examples are a column-generation based approach 
for robust knapsack problems and shortest path problems with uncertain demand \cite{BAH-ESA11}, 
an approach using Bender's decomposition for railway rolling stock planning \cite{Kroon10},
and the idea of replacing the recovery algorithm by a metric \cite{Tapas2011-Goerigk-Sch,GoeSchoe13a,marcdiss}.

\medskip

\emph{Example 3:} Robust passenger information systems.

The following example shows the application of the algorithm engineering cycle on
a specific application, namely constructing robust timetable information systems.
Suppose that a passenger wants to travel from an origin to some destination by public transportation. The
passenger can use a timetable information system which will provide routes with small traveling time. However,
since delays are a matter of fact in public transportation, a robust route would be more valuable than just
having a shortest route. In \cite{rotti} this problem was considered for strictly robust routes: The model was
set up, analyzed (showing that it is NP-complete), and an algorithm for its solution was designed. 
The experimental evaluation on real-world data showed that strictly robust routes are useless in practice: their
traveling time is much too long. Based on these experiments, light robust passenger information system was
considered. The light robust model was designed and analyzed; algorithms based on the strictly robust
procedures could be developed. The experiments showed that this model is much better applicable in practice.
However, the model was still not satisfactory, since it was assumed that a passenger stays on his/her route whatever happens.
This drawback motivated to start the algorithm engineering cycle again in \cite{rotti2} 
where now recoverable robust
timetables are investigated.
\bigskip

Considering the cycle of
design, analysis, implementation, and experiments, we were also able to
identify pointers for further research. We summarize our results by
reproducing the most significant messages:

\begin{enumerate}
\item {\it Robust optimization is application-driven.} From the beginning, robust optimization was intended as an optimization approach which generates solutions that perform well in a realistic environment. As such, it is highly appealing to practitioners, who would rather sacrifice some nominal solution quality if the solution stability can be increased. 
\item {\it Robust optimization needs a unified classification scheme.} While the strong connection to applications is a beneficial driver of research, it also carries problems. One striking observation is a lack of unification in robust optimization. This begins with simple nomenclature: The names for strict robustness, or the uncertainty set considered by Bertsimas and Sim are plenty. It extends to the frequent re-development of algorithmic ideas (as iterative scenario generation), and the reinvention of robustness concepts from scratch for specific applications. This lack of organization is in fact unscientific, and endangers the successful perpetuation of research. As related problems, some journals don't even offer ``robust optimization'' as a subject classification (even though publishing papers on robust optimization); solutions generated by some fashion that is somehow related to uncertainty call themselves ``robust''; and students that are new to the field have a hard time 
to identify the state-of-the-art.
\item {\it Performance guarantees are not sufficiently researched in robust optimization.} Also this point can be regarded as related to robust optimization being application-driven and non-unified. Performance guarantees are of special importance when comparing algorithms; hence, with a lack of comparison, there also comes a lack of performance guarantees. This includes the comparison of robust optimization concepts, of robust optimization algorithms, and even the general evaluation of a robust solution compared to a non-robust solution.
\item {\it There is no robust optimization library available with specifically designed algorithms other than reformulation approaches.}
While libraries for robust optimization exist, they concentrate on the modeling aspects of uncertainty, and less on different algorithmic approaches. Having such a library available would 
prove 
extremely helpful not only for practitioners, but also for researches that develop new algorithms and try to compare them to the state-of-the-art.
\item {\it There are too few comparative studies in robust optimization.} All the above points culminate in the lack of comparative studies; however, we argue that here also lies a chance to tackle these problems. This paper is a humble step to motivate
such research, and we hope for many more publications to come.
\end{enumerate}

\bibliographystyle{alpha}

\newcommand{\etalchar}[1]{$^{#1}$}

\end{document}